\documentclass{amsart}
\usepackage{amssymb,amsfonts,amsmath}
\usepackage[all]{xy}
\usepackage[shortlabels]{enumitem}
\usepackage{hyperref, color}		

\newtheorem{theorem}{Theorem}[section]
\newtheorem{proposition}{Proposition}[section]

\theoremstyle{definition}
\newtheorem{definition}{Definition}[section]
\newtheorem{lemma}{Lemma}[section]

\theoremstyle{remark}
\newtheorem{remark}{Remark}
\theoremstyle{case}

\theoremstyle{claim}
\newtheorem{claim}{Claim}[section]

\theoremstyle{theoremA}

\theoremstyle{theoremB}
  
\newcommand{\real}{\mathbb{R}}

\newcommand{\Sp}{\mathbb{S}}

\newcommand{\om}{\omega}

\newcommand{\ga}{\gamma}

\newcommand{\ti}{\tilde}

\newcommand{\lan}{\left\langle}
\newcommand{\ran}{\right\rangle}

\newcommand{\grad}{\mathrm{grad}}

\newcommand{\mk}{\mathfrak}

\newcommand{\ak}{\check{\alpha}}
\newcommand{\nk}{\check{\nabla}}
\newcommand{\setap}{\rightarrow}
\begin{document}

\title[Submanifolds into Rotational Hypersurfaces]{A Bonnet Theorem for Submanifolds into Rotational Hypersurfaces}
\thanks{The first author was partially supported by CAPES-Brazil}

\author{C. do Rei Filho}

\author{F. Vit\'orio}

\date{}

\maketitle
\thispagestyle{empty}
\begin{abstract} In this work, we prove a version of the fundamental theorem of submanifolds to target manifolds with warped structure. 
\end{abstract}
\section{Introduction}

The classical Bonnet's theorem establishes necessary and
sufficient conditions for the existence of an isometric immersion
of a simply connected manifold $M^k$ in an Euclidean space
$\mathbb{R}^{d+l}$ with prescribed second fundamental form. There
are several proofs for this beautiful achievement in Differential
Geometry. One could read for example the original proof for
$\real^3$ in \cite{bonnet}, or more general versions in
\cite{spivak}, \cite{petersen} or \cite{marcos}. It seems that the
underlying idea of the Bonnet's proof is that the metric and the
second fundamental form of the immersed manifold must obey some
compatibility equations. 

The recent development of the study of constant mean curvature
surfaces in the three dimensional geometries with $3$ and $4$
dimensional isometries groups exploits the resemblances with the
known results in the theory on space forms. Benoit Daniel,
\cite{dan} and \cite{daniel2}, had proved versions of Bonnet's
theorem in $\mathbb{S}^2\times \real$, $\mathbb{H}^2\times \real$,
three dimensional Berger spheres and three dimensional Heisenberg
groups. His method is to add to the Gauss and Codazzi equations
some conditions about the tangential and normal components of a
``vertical'' vector field and to verify directly that these
conditions fit the compatibility equations for integrate a certain
distribution. As far the authors know, \cite{dan} is the first paper that deals with the problem which the classical structure equations are not sufficients to give the immersions. In fact, the tangent and normal projections yield  additional  conditions for the immersion. Along the recent year some  works  \cite{cx}, \cite{{Kowal}} and \cite{ltv} gave generalizations for Daniel's work in some directions.
The aim of this work is give another piece in this mosaic solving the problem for a class of non-homogeneous target manifolds. Our approach, like \cite{ltv}, is to reduce the problem to the case of flat Euclidean
target spaces by using the canonical isometric embedding of a revolution hypersurface into a flat Euclidean space. We point out that  in \cite{cx} the authors prove a version of a Bonnet's theorem for warped products provided that the basis is a flat space, i.e., $\real\times_\eta\real^n$.

Our main result uses a fiber bundle terminology, see Definition \ref{compatibility} and Theorem \ref{main} for the precise statement, but it can be read roughly as following:
\vspace{0.4cm}

\begin{minipage}[b]{11.5cm}
  {\it
Let $\big(M^k,g \big)$ be a simply connected
Riemannian manifold of dimension $k$ and let $\bar M$ be a revolution hypersurface.  There exist 
compatibility equations for $\bar M$  that are necessary and sufficients for  existence of  an isometric immersion $
\mathbf{x}:M^k \to
\bar M^n$. Furthermore, the isometric immersion induces a vector bundle isomorphism.
}
\end{minipage}

\section{Preliminaries.}\label{prelim}
Let $\mathbb{E}^n$ be the $\real^n$ with metric $(1,...,1,\epsilon),$ where $\epsilon \in \{-1,1\}.$ 
Let $\ga: I \to \mathbb{E}^2$ be a curve parametrized by the arc lenght,  \[ \ga(t)=\big(f(t), h(t)\big),\]
where $f,h:I \to \real$ are smooth function satisfying $f'(t)^2+\epsilon h'(t)^2=1$, for all $t\in I$. For our purposes, we will assume that $f$ is positive and $h'(t)>0$.
  Let $\Phi: I \times \Sp^n \to \mathbb{E}^{n+2}$ be the rotational hypersurface given by
\[
\Phi(t,\om)=\big( f(t) \om, h(t)\big),
\]
where $\om \in \Sp^n$ and we are considering via the canonical embeddings $\Sp^n\subset \real^{n+1}= \real^{n+1}\times\{0\}\subset \mathbb{E}^{n+2}$.  On the cylinder $I\times \Sp^n$ we define the Riemannian metric $ ds^2= \Phi^*(\lan\,,\ran)$, where $\lan\,,\ran$ is the standard metric of $\mathbb{E}^{n+2}$. It is simple to see that $\bar M^{n+1}=\big(I\times \Sp^n, ds^2\big)$ is a warped product manifold with warped metric $ds^2= dt^2+f(t)^2d\sigma^2$ where $d\sigma^2$ is the standard metric of $\Sp^n$. Note that, $\bar M$ has a distinguished unitary vector field $\partial_t= (f'(t)\om, h'(t))$. 

Notice that the unitary normal field of $\Phi$ is given by $N_t=(h'(t)\omega,-\epsilon f'(t))$. Furthermore, if we consider the map $\ti{\Psi}:I\times \Sp^n \to \mathbb{E}^{n+2}$ defined by  $$\ti{\Psi}(t,\omega)=\Phi(t,\omega)-f(t).f'(t)\partial_t - \epsilon f(t).h'(t)N_t=(0,...,0,h(t)).$$
Then, $\ti \Psi$ does not depends on $\omega$ and, hence, the curve  $\ti \sigma:I \to \mathbb{E}^{n+2}$, given by  $t\mapsto \ti{\Psi}(t,\omega) $ is well defined. See also that,
\begin{equation*} 
\ti \sigma'(t)=(0,...,0,h'(t))=\epsilon h'(t)\ [ h'(t)\partial_t-f'(t)N_t].
\end{equation*}
Conversely, let $\bar M^{n+1}=\big(I\times \Sp^n, ds^2\big)$ be a warped product manifold with warped metric $ds^2= dt^2+f(t)^2d\sigma^2$ where $d\sigma^2$ is the standard metric of $\Sp^n.$ If $(1-f'(t)^2)\epsilon >0, \ \forall \, t \in I$ then is well defined a function $h:I \to \real$ by the expressions $f'(t)^2+\epsilon h'(t)^2=1$ and$\ h'(t)>0$, for all $t\in I$. Therefore  $\Phi: \bar M^{n+1} \to \mathbb{E}^{n+2}$ defined by
\begin{equation}\label{isometric}
\Phi(t,\om)=\big( f(t) \om, h(t)\big),
\end{equation}
is an isometric immersion. \\

Let $\bar M=\big(I\times \Sp^n, ds^2\big)$ be a warped product manifold with  metric $ds^2= dt^2+f(t)^2d\sigma^2$. Let $\bar{\nabla}$  be the  Levi-Civita connection on   $\bar M$. The covariant derivative of the vector field $\partial_t$ with respect to any tangent vector satisfies
\begin{equation}
\bar{\nabla}_u \partial_t = \frac{f'(t)}{f(t)}(u-\langle u,\partial_t \rangle \partial_t), \,\,\, \forall \, u \in T\bar{M}.
\end{equation}  Thus, in particular, the  orbits of the vector field $\partial_t$ are geodesics and the vector field $V=f(t) \partial_t$ is closed conformal with conformallity factor $f'(t)$ i.e., 
\begin{eqnarray}\label{c.c.f.}
\bar\nabla_u V=f'(t)u, \,\, \forall \, u \in T\bar{M}.
\end{eqnarray}
It is clear that $\bar{M}$ is foliated by spheres $\Sigma_t=\Sp^n(f(t))$ for each $t\in I$, $\partial_t$ is the unit normal vector field for each leaf $\Sigma_t$ which is umbilical, moreover  $\Sigma_t$ has mean curvature $-\frac{f'(t)}{f(t)}$ with sectional curvature equals to $1/f(t)^2$.

Let  $\bar{R}$ be the curvature tensor on
$\bar M$. We observe that
\begin{eqnarray}
\bar{R}(u,v)\partial_t = \frac{f''(t)}{f(t)}\big(\langle u,\partial_t \rangle v - \langle v,\partial_t \rangle u\big), \,\,\, \forall \, u,v \in T\bar{M}.
\end{eqnarray}

In this way, we can compute the curvature tensor $\bar{R}$ as follows

 \begin{equation} \label{curv-tensor}
 \begin{array}{rcl}
\bar{R}(u,v)w &=&  \Big(\frac{1-f'(t)^2}{f(t)^2} + \frac{f''(t)}{f(t)} \Big) \Big[ \langle w,\partial_t  \rangle\big( v,u ,\partial_t \big) + \langle \big(v,w,u \big),\partial_t \rangle  \partial_t \Big] \\ 
                           &&+ \frac{1-f'(t)^2}{f(t)^2} \big( u,v,w\big) 
\end{array}                           
\end{equation}
where $\big( u,v,w\big)=\langle v,w \rangle u - \langle u,w \rangle v$, for all $u,v,w \in T\bar{M}$.

It is worthwhile to mention that the shape operator $A$ can be expressed as $$Av=-\frac{h'(t)}{f(t)}v+\Big(\frac{h'(t)}{f(t)}+\epsilon\frac{f''(t)}{h'(t)}\Big)\langle v,\partial_t \rangle \partial_t \quad \forall \,\,v\in T\bar M.$$

\section{Necessary conditions for an isometric immersion}
Let $\mathbf{x}:M^k \to \bar M$ be an isometric immersion of a $k-$dimensional Riemannian manifold $M^k$ into the warped product $\bar M$. Let us denote $T$ the canonical projection of $\partial_t$ into $TM,$ such that
\begin{equation}\label{troc.xino}
\partial_t=T + \varrho,
\end{equation} 
where $\varrho $ is a section of $TM^{\perp}.$

\begin{proposition}\label{prop.princ.1}Let $\mathbf{x}:M^k \to \bar M$ be an isometric immersion of a $k-$dimensional Riemannian manifold $M^k$ into the warped product $\bar M$.
For all $u \in TM,$ we have
\begin{equation}\label{bola}
\alpha(u,T) + \nabla^{\perp}_u \varrho = -\frac{f'(t)}{f(t)}\langle u, T \rangle \varrho
\end{equation}
\begin{equation}\label{lopo}
\nabla_u T - A_{\varrho}u= \frac{f'(t)}{f(t)}(u-\langle u,T \rangle T).
\end{equation}
where $\nabla,$ $\nabla^{\perp},$ $\alpha$ and $A$ denote the Levi-Civita conection in $M,$ the induced normal conection in $M,$ the second fundamental form of the immersion $\mathbf{x} (M)$  and the Weingarten endomorphism associated to $\alpha,$ respectively.
\end{proposition}
\noindent{\bf Proof:}
Derive $\partial_t=T + \varrho $ with respect to $u$ and use the Gauss and Weingarten formulae. We conclude the result taking the tangent and normal components of that derivative.
\begin{flushright}$\Box$
\end{flushright}

The expression of the tensor curvature (\ref{curv-tensor}) allow us write the Gauss, Codazzi and Ricci equations in the following way:
\begin{proposition}\label{prop.princ.2}
Let $\mathbf{x}:M^k \to \bar M$ be an isometric immersion of a $k-$dimensional Riemannian manifold $M^k$ into the warped product $\bar M$. 
 Let us consider $u,v,w,z,w \in TM$ and $\xi, \eta \in TM ^\perp$. Then, Gauss, Codazzi and Ricci equations, for the isometric immersion $\mathbf{x}$, are, respectively:

\begin{eqnarray}
\langle R(u,v)w,z \rangle \! &=& \! \lambda(t) \big( u,v,w,z\big) \nonumber \\ 
                          & &+\,\, \mu(t) \big[  \lan w, T\ran \big(v,u,T,z\big)+\big( v,w,u,T\big) \lan T,z \ran \big]\\
                          & &+\,\, \langle \alpha(u,z), \alpha(v,w) \rangle - \langle \alpha(u,w),\alpha(v,z) \rangle \nonumber
\end{eqnarray}

\begin{eqnarray}
(\nabla_u ^{\perp} \alpha)(v,w) - (\nabla_v ^{\perp} \alpha)(u,w) = \mu  \big( v,w,u,T\big) \varrho
\end{eqnarray}

\begin{equation}
\langle R^\perp (u,v)\xi,\eta \rangle = \langle [A_{\xi},A_{\eta}]u,v \rangle,
\end{equation}

where \begin{equation*}
\lambda(t)= \frac{1-f'(t)^2}{f(t)^2} \quad \textrm{and} \quad \mu(t) = \frac{1-f'(t)^2}{f(t)^2} + \frac{f''(t)}{f(t)} .
\end{equation*}
\end{proposition}
\noindent{\bf Proof:} It is a standard computation and follows directly from  (\ref{curv-tensor}).
\begin{flushright}$\Box$
\end{flushright}

To finish this section we want point out an intrisical characterisation of the vector field $T$. This will be given by the following
\begin{proposition}\label{gradient} Let $\mathbf{x}:M^k \to \bar M$ be an isometric immersion of a $k-$dimensional Riemannian manifold $M^k$ into the warped product $\bar M$.  
 There exists a function $\mk h: M \to \real$ such that $T=\grad_M \mk h$.
\end{proposition}
 
 \section{Establishing the sufficient conditions}

Let  $(M^k, \langle\,,\,\rangle)$ be a Riemannian manifold and let us denote  $\nabla$ its Levi-Civita connection.   Let $E$ be a Riemannian vector fiber bundle on $M^k$ with rank $n+1-k$ and let us denote  $\nabla'$ its compatible connection. Let $\mk h \in C^{\infty}(M),\, \varrho \in \Gamma(E)$ and $\alpha '$  be a smooth function on $M$, a section of the vector fiber bundle $E$ and a symmetric section of the homomorphisms fiber bundle $\textrm{Hom\,}(TM \times TM, E)$, respectively. Let us define for each local section  $\xi\in \Gamma(E)$ the map $A'_{\xi}:TM \to TM$ by 
\begin{equation}
\langle A'_{\xi}u,v \rangle = \langle \alpha'(u,v),\xi \rangle, 
\end{equation}
for all $u,v \in TM$.

In virtue of the proposition \ref{gradient} we can rephrased  the necessary conditions obtained in the  propositions \ref{prop.princ.1} and \ref{prop.princ.2} using this abstract framework of fiber bundles.

\vspace{0.1cm}

\begin{definition}\label{compatibility}
We say that the data $(M^k,\langle,\rangle, \nabla',\alpha',\varrho,\mk h)$ satisfies the compatibility equations for $\bar{M}$ if 

\begin{equation}\label{eq.comp.1}
{\arrowvert {T} \arrowvert}^2 + {\arrowvert \varrho \arrowvert}^2 =1, \,\,\,\,\, {T=\grad_M \mk h}
\end{equation}
and for all $u,v,z,w \in \chi(M)$ and $\xi, \eta \in \Gamma( E)$ the following equations hold:
\begin{equation}\label{eq.comp.2}
\alpha '(v,T) + \nabla'_v \varrho =  -\frac{f'(\mk h)}{f(\mk h)}\langle v, T \rangle \varrho
\end{equation}

\begin{equation}\label{eq.comp.3}
\nabla_v T - A'_{\varrho}v= \frac{f'(\mk h)}{f(\mk h)}\big( v-\langle v, T \rangle X \big)
\end{equation}
\begin{eqnarray}\label{eq.comp.4}
\langle R(u,v)w,z \rangle \! &=& \! \lambda(\mk h) \big( u,v,w,z\big) \nonumber \\   & &+\,\, \mu(\mk h) \big[  \lan w, T\ran \big(v,u,T,z\big)+\big( v,w,u,T\big) \lan T,z \ran \big]\\
                          & &+\,\, \langle \alpha(u,z), \alpha(v,w) \rangle - \langle \alpha(u,w),\alpha(v,z) \rangle \nonumber
\end{eqnarray}

\begin{eqnarray}
(\nabla_u ^{\perp} \alpha)(v,w) - (\nabla_v ^{\perp} \alpha)(u,w) = \mu(\mk h) \big( v,w,u,T\big) \varrho
\end{eqnarray}

\begin{equation}
\langle R^\perp (u,v)\xi,\eta \rangle = \langle [A_{\xi},A_{\eta}]u,v \rangle,
\end{equation}
where \begin{equation*}
\lambda(\mk h)= \frac{1-f'(\mk h)^2}{f(\mk h)^2} \quad \textrm{and} \quad \mu(\mk h) = \frac{1-f'(\mk h)^2}{f(\mk h)^2} + \frac{f''(\mk h)}{f(\mk h)} .
\end{equation*}
\end{definition}

\section{Proof of the Fundamental Theorem }

Let us consider the fiber bundle obtained by the  Whitney sum of the tangent fiber bundle $TM$ with  the fiber bundle $E$ as in the definition \ref{compatibility}, $\ti{E}= TM\oplus_w E$, endowed with the product metric and compatible connection
\begin{eqnarray*}
\nabla''_vu &=& \nabla_vu + \alpha'(v,u) \quad u,v\in TM ,\\
\nabla''_v\xi &=& -A'_{\xi}v + \nabla'_v\xi \quad v \in TM \,\,and \,\, \xi \in \Gamma(E).
\end{eqnarray*}

It is easy see that, the section $X=T+\varrho \in \Gamma (\ti{E})$ satisfies $\arrowvert X \arrowvert = 1$ and
\begin{equation}
\nabla''_vX=\frac{f'(\mk h)}{f(\mk h)}(v-\langle v,T \rangle X), \quad \forall v\in TM.
\end{equation}
In particular,
\begin{equation*}
\nabla''_vf(\mk h)X=f'(\mk h)v, \quad \forall v\in TM.
\end{equation*}

Moreover, if $\langle v,T \rangle = \langle v, X \rangle = 0,$ then  $v(f(\mk h))=v(h(\mk h))= v(\mk h)=0,$ since $v(\mk h)=\langle \grad_M \mk h,v\rangle = \langle T,v\rangle. $ \\

Now, assume that $(1-f'(\mk h)^2)\epsilon >0$. Thus, the  function $h$ given by the expressions $f'(\mk h)^2+\epsilon h'(\mk h)^2=1$ and$\ h'(\mk h)>0$ is well defined, up to a constant.
Let $\check{E}={E}\oplus_w\langle \zeta \rangle$ be the Semi-Riemannian fiber bundle obtained by summing the Semi-Riemannian line fiber bundle  $\langle \zeta \rangle$ to ${E}$. On $\check{E}$, we define
\begin{eqnarray*}
\nk \!\!\!\! &:& \!\!\!\! TM\times \check{E} \setap \check{E} \\
\ak \!\!\!\! &:& \!\!\!\! TM\times TM \setap \check{E}
\end{eqnarray*} 
putting,
\begin{eqnarray}
\ak (u,v) &=& \alpha' (u,v) +\epsilon(- \ti \lambda (\mk h) \langle u,v \rangle
 + \ti \mu (\mk h) \langle u,X\rangle \langle v,X \rangle)\zeta \\
\nk_v\phi &=& \nabla'_v(\phi)_E +\epsilon \ti \mu (\mk h) \langle v,X \rangle \Big( \langle \phi,X \rangle \zeta - \langle \phi,\zeta \rangle \varrho\Big)  +\epsilon v(\langle \phi,\zeta\rangle)\zeta
\end{eqnarray}
where $u,v \in TM$,  $\phi \in \Gamma(\check{E})$,  $(\phi)_E$ is the canonical projection on $E$ and $\ti \lambda, \ti \mu$ are defined by
\begin{equation}
\ti \lambda(\mk h)= \frac{h'(\mk h)}{f(\mk h)} \quad \textrm{and} \quad \ti \mu (\mk h)= \frac{h'(\mk h)}{f(\mk h)} + \epsilon\frac{ f''(\mk h)}{h'(\mk h)} 
\end{equation}

\begin{remark}
Note that $\ti \lambda (\mk h)^2=\epsilon\lambda(\mk h)$ and $\ti \lambda(\mk h).\ti \mu (\mk h)= \epsilon\mu(\mk h).$
\end{remark}

Under the notations and definitions above is straightforward conclude the following 
\begin{lemma}
Assume that the data $(M^k,\langle,\rangle, \nabla',\alpha',\varrho,\mk h)$ satisfy the compatibility equations for­ $\bar{M}$. If $(1-f'(\mk h)^2)\epsilon >0$, then the data $(M,g,\nk,\ak)$ satisfy the compatibility equations for $\mathbb{E}^{n+2}$.
\end{lemma}

Thus, using the fundamental theorem of submanifolds, there exists an isometric immersion $g:M^k \setap \mathbb{E}^{n+2}$ and a fiber bunder isometry $\check{g}:\check{E} \setap TM^\bot$ along $g$, such that 
\begin{eqnarray}
\ti \alpha = \check{g} \check{\alpha}
\end{eqnarray}
\begin{eqnarray}
\ti{\nabla}^{\bot} \check{g}  =\check{g} \check{\nabla}  \nonumber
\end{eqnarray}
where $\ti{\nabla}^{\perp}$ and $\ti{\alpha}$ are the normal connection and second fundamental form of  $g(M) \subset \mathbb{E}^{n+2},$ respectively.

Denoting by $D$ the covariant derivative of $\mathbb{E}^{n+2}$, it is a simple computation, identifying $\check{g}(X)$ with $X$ and $\check{g}(\zeta)$ with $\zeta$ show

\begin{equation} \label{equ.deriv}
\begin{array}{rcl}
D_v X &=& \frac{f'(\mk h)}{f(\mk h)}(v - \langle v,T \rangle X) + \epsilon
\langle v,T \rangle(\ti \mu - \ti \lambda)\zeta, \quad \forall \, v \in TM. \\ \\
D_v \zeta &=& \ti \lambda v - \ti \mu \langle v, T \rangle X, \quad \forall \, v \in TM. \nonumber
\end{array}
\end{equation}

\begin{claim} \label{prop.frob} 
For all $u,v\in TM,$ with $\langle u,T \rangle = \langle v,T \rangle =0,$ we have that $\langle D_vu,X\rangle = \langle D_uv, X \rangle.$ In particular, $\langle [u,v],T \rangle =0. $
\end{claim}

 Now, assume further that  $\mk h:M\setap \real$ is a smooth submersion, i.e., $T=\grad_M \mk h \neq 0$ for all point in $M$.


First of all, note that the claim
\ref{prop.frob} implies that the distribution 
$$p\in M \longmapsto \mathfrak{D} (p) =\{v\in T_pM; \langle v, T \rangle =0\}$$
is  involutive, hence totally integrable.
 
Therefore, $M$ admits a codimension one foliation  $\mathfrak{F}(T)$ oriented by the unitary vector field  $T/|T|$.  Note also that the foliation $\mathfrak{F}(T)$ is determined by the submersion $\mk h:M \setap \real$ as level sets, on the other words, $\mathfrak{F}(T)=\{\mk h^{-1}(s)\subset M;s\in \mk h(M)\}$.\\

Let $\Psi:M\setap \mathbb{E}^{n+2}$ be the smooth map defined by
 $$\Psi(p)=g(p)-f(\mk h(p)).f'(\mk h(p))X|_p - \epsilon f(\mk h(p)).h'(\mk h(p))\zeta|_p.$$

\begin{claim}
The map $\Psi$ is constant along each connected  leaf of the foliation $\mathfrak{F}(T).$
\end{claim}
\noindent {\bf Proof:}  In order to compute a derivative of $\Psi$, note that for $v \in \mathfrak{D}$ the equations (\ref{equ.deriv}) can be rewritten as
\begin{equation*}
\begin{array}{rcl}
D_v X &=& \frac{f'(\mk h)}{f(\mk h)}v,  \\ \\
D_v \zeta &=& \ti \lambda v=\frac{h'(\mk h)}{f(\mk h)} v.
\end{array}
\end{equation*}
Hence, 
\[
D_v \Psi= v- f(\mk h(p)).f'(\mk h(p))D_v X - \epsilon f(\mk h(p)).h'(\mk h(p))D_v\zeta =v- f'^2 v - \epsilon h'^2 v=0,
\]
since that $f'^2+\epsilon h'^2=1$. 
\begin{flushright}$\Box$
\end{flushright}

Now, the connectedness of $M$ and the fact that $\mk h:M \setap \real$ is a smooth submersion imply that $\mk h(M)\subset \real$ is an open interval. Moreover, since the map $\Psi$ is constant along each level set of $\mk h$, there exists an unique smooth map $\sigma:\mk h(M)\subset \real \setap \mathbb{E}^{n+2}$ such that $\Psi=\sigma\circ \mk h.$ 
\begin{center}
\hspace{0.5cm}\xymatrix{
M  \ar[r]^{\Psi}
\ar[d]^{\mk h} & \mathbb{E}^{n+2} \\
\mk h(M) \ar[ur]_{\sigma} }
\end{center}
 For each $s\in \mk h(M),$ the equation 
$$\langle \sigma(s)-g(p),\sigma(s)-g(p)\rangle=f(s)^2,\,\, \forall \, p\in \mk h^{-1}(s),$$ 
shown that image of each leaf $\mk h^{-1}(s)$ is contained in a $(n+1)-$dimensional pseudosphere of $\mathbb{E}^{n+2}$ centered at $\sigma(s)$ and ratio $f(s)$, i.e., $g(\mk h^{-1}(s)) \subset \mathbb{S}^{n+1}_{f(s)}(\sigma(s))\subset \mathbb{E}^{n+2}$. We call $\sigma:\mk h(M) \setap \mathbb{E}^{n+2}$ the {\em curve of centers} of $M$.


\begin{claim}\label{cenret}
$\sigma(\mk h(M))\subset \mathbb{E}^{n+2}$  is a straight line. 
\end{claim}
\noindent{\bf Proof:} In order to prove this claim, we will prove that the curvature of the curve of centers is identically zero. Since that $\Psi=\sigma\circ \mk h$ and $T(\mk h)=|T|^2$, we obtain $\sigma'(\mk h)=|T|^{-2}T(\Psi).$  Thus, computing $T(\Psi),$ we have $$\sigma'(\mk h)= \epsilon h'(\mk h)\Big( h'(\mk h)X-f'(\mk h)\zeta\Big ).$$
In the same way, $\sigma''(\mk h)=|T|^{-2}T(\Psi')$, where $\Psi'=\sigma'\circ \mk h$. Thus, computing $T(\Psi')$ we have
$$\sigma''(\mk h)= \epsilon h''(\mk h)\Big( h'(\mk h)X-f'(\mk h)\zeta\Big ).$$
Thus, 
\begin{equation} \label{zerocurvature}\sigma''(\mk h)= \frac{h''(\mk h)}{h'(\mk h)} \sigma'(\mk h).
\end{equation}
Therefore, (\ref{zerocurvature}) implies that the curvature of $\sigma(\mk h)$ is zero, since that the velocity of the curve $\sigma(\mk h)$ is $|\sigma'(\mk h)|=h'(\mk h)$. \begin{flushright}$\Box$ \end{flushright}


Note that $\langle \sigma',\sigma' \rangle = \epsilon h'^2.$ As a direct consequence of the claim \ref{cenret} we have
 $$H=\{x\in \mathbb{E}^{n+2} ; \langle x,\sigma'(s)  \rangle =0\}$$  is a Riemannian hyperplane of  $\mathbb{E}^{n+2}$ which does not depends of the parameter  $s$.
Thus the equation $$\langle \sigma'(s),\Psi(p)-g(p)\rangle = 0, \, \forall \, p\in \mk h^{-1}(s), $$ tell us that for each  $s\in \mk h(M),$ $$g(\mk h^{-1}(s)) \subset \big(\mathbb{S}^{n+1}_{f(s)}(\sigma(s))\, \cap \,(\sigma(s)+ H) \big) .$$
\medskip

Let us take $s_0\in \mk h(M)$ and let $\tau:\mathbb{E}^{n+2}\setap \mathbb{E}^{n+2}$ be the rigid motion of $\mathbb{E}^{n+2}$ such that $\tau(\sigma(\mk h(M)))$ is contained in the axis ${\Huge O}x_{n+2}$, $(0,...,0,h(s_0))=\tau(\sigma(s_0))$ and the velocity vector $\tau(\sigma')$ is pointing in the same orientation of axis  ${\Huge O}x_{n+2}$. Note that such isometry satisfies $$\tau(g(\mk h^{-1}(s_0))) \subset \big(\mathbb{S}^{n+1}_{f(s)}((0,...,0,h(s_0)))\, \cap \,((0,...,0,h(s_0))+\{x\in \mathbb{E}^{n+2};x_{n+2}=0\}) \big).$$
Thus, for construction, the curves  $s\in \mk h(M)\setap (0,...,0,h(s))$ and $s\in \mk h(M)\setap \tau(\sigma(s))$ coincide in $s_0$, their derivatives coincide for all points, ($\tau(\sigma'(s))$ and $h'(s)$ are pointing in the same direction and $|\tau(\sigma'(s))|=h'(s)$), hence $(0,...,0,h(s))=\tau(\sigma(s)), \forall \, s\in \mk h(M)\subset I.$ Therefore, we conclude that­ $\tau(g(M))\subset \Phi(\bar M^{n+1}),$ where $\Phi$ is the isometric gimmersion given in (\ref{isometric}). 

Now, we will summarize the informations that we have obtained, up to an isometry of $\mathbb{E}^{n+2}$,
\begin{enumerate}
\item[a)] There exists an isometric immersion  $g:M^k \setap \mathbb{E}^{n+2}$ and a fiber bundle isometry $\check{g}:\check{E} \setap TM^\bot$ along $g$, such that $\ti \alpha = \check{g} \check{\alpha}$ and $\ti{\nabla}^{\bot} \check{g}  =\check{g} \check{\nabla},$ where $\ti{\nabla}^{\perp}$ and $\ti{\alpha}$ are the normal connection and second fundamental form of  $g(M) \subset \mathbb{E}^{n+2},$ respectively.\\

\item[b)]  The curve $\sigma(s)=g(p)-f(s).f'(s)X|_p -\epsilon f(s).h'(s)\zeta|_p$   parametrize an open interval of the axis ${\Huge O}x_{n+2}$. More specifically,  $\sigma(s)=(0,...,0,h(s)).$  We also have $\sigma'(s)=\epsilon h'(s)\Big(h'(s)X - f'(s)\zeta\Big )=(0,...,0,h'(s)).$\\
\item[c)] We have, for each  $s\in \mk h(M)$, $$g(\mk h^{-1}(s)) \subset \big(\mathbb{S}^{n+1}_{f(s)}(\sigma(s))\, \cap \,(\sigma(s)+\{x\in \mathbb{E}^{n+2};x_{n+2}=0\}) \big).$$ In particular, $g(M)\subset \Phi(\bar M^{n+1}).$\\
\item[d)] Furthermore, $(0,...,0,h(s))=\Phi(s,\omega)-f(s).f'(s)\partial_t - f(s).h'(s)N_t $ and $(0,...,0,h'(s))=h'(s)^2\partial_t-f'(s)h'(s)N_t.$\\
\item[e)] The items (b) and (d) provide that $X=\partial_t|_{g(M)}$  and $\zeta=N_t|_{g(M)}$. Indeed, given $p \in M$, take  $\omega \in \mathbb{S}^n$ such that  $g(p)=\Phi(s,\omega)$. Thus, the items (b) e (d) above, imply in the following system:
\begin{equation}
\left\{ \begin{array}{l} \label{sys1}
-f(s).f'(s)X - f(s).h'(s)\zeta = -f(s).f'(s)\partial_t - f(s).h'(s)N_t \\ \\

h'(s)^2X - f'(s)h'(s)\zeta = h'(s)^2\partial_t - f'(s)h'(s)N_t.
\end{array}\right.
\end{equation}
Since that  $f(s)\neq 0$ and $h'(s)\neq 0, \,\, \forall\, s\in I$, the system of equations (\ref{sys1}) is equivalent to
\begin{equation} \label{sys2}
\left\{ \begin{array}{l}
f'(s)(X-\partial_t) = -h'(s)(\zeta - N_t)\\\\
h'(s)(X-\partial_t) = f'(s)(\zeta - N_t)
\end{array}\right.
\end{equation}
Note that the system of equations (\ref{sys2}) implies that $X=\partial_t$ and $\zeta=N_t$ at $p$. As the choice of $p\in M$ was arbitrary, the result follows.
\end{enumerate}

Therefore, there exists an isometric immersion  $\mathbf{ x}:M^k\setap \bar M,$ defined by $g=\Phi \circ \mathbf{ x}$ and a fiber bundle isometry $\mathbf{\ti x}: E \setap TM^\bot$ along $\mathbf{ x},$  defined by $\mathbf{\ti x}=\check{g}|_E,$ such that $\alpha=\mathbf{\ti x} \alpha'$ and $\nabla^{\bot} \mathbf{\ti x}=\mathbf{\ti x}\nabla',$
where $\nabla^{\perp}$ and $\alpha$ are the normal connection and the second fundamental for of $\mathbf{ x}(M) \subset \bar M,$ respectively.
Moreover, $$\partial_t=\mathbf{x}_*(T)  + \mathbf{\ti x} (\varrho).$$ 



This allow us to conclude the following

\begin{theorem}\label{main}
Let $\big(M^k,\langle,\rangle \big)$ be a $k-$dimensional simply connected Riemannian manifold. Assume that the data $(M^k,\langle,\rangle, \nabla',\alpha',\varrho,\mk h)$ satisfy  the compatibility equations for $\bar M$, as in the definition \ref{compatibility}, assume also that  $\mk h$ is a smooth submersion and that $(1-f'(\mk h)^2)\epsilon >0$.  Then, there exists an isometric immersion $\mathbf{x}:M^k \to
\bar M^{n+1}$  and  a fiber bundle isometry  $ \mathbf{\ti x}: E
\rightarrow TM^\bot$ along  $\mathbf{x}$,  such that  $$\alpha=\mathbf{\ti x}\alpha' \quad {e} \quad \nabla^{\bot}\mathbf{\ti x}= \mathbf{\ti x}\nabla' \quad \textrm{and} \quad \partial_t=\mathbf{x}_*(\grad_M \mk h) + \mathbf{\ti x}(\varrho)$$ 
where $\nabla^{\bot}$ and $\alpha$ are the normal connection and the second fundamental form of $\mathbf{x}(M) \subset \bar M$, respectively.
\end{theorem}

\vspace{0.7cm}
\begin{small}

\begin{tabular}{l}
Carlos do Rei Filho\\
Universidade Federal de S\~ao Carlos\\
Departamento de Matem\'{a}tica\\
13565-905 S\~ao Carlos-SP\\
Brazil\\
\verb+carlosfilho@dm.ufscar.br+\\
\end{tabular}\

\vspace{1cm}

\begin{tabular}{l}
Feliciano Vit\'orio\\
Universidade Federal de Alagoas\\
Instituto de Matem\'{a}tica\\
57072-900 Macei\'o-AL\\
Brazil\\
\verb+feliciano@pos.mat.ufal.br+\\
\end{tabular}\

\end{small}

\end{document}